\documentclass{article}
\usepackage{array,amsmath}
\numberwithin{equation}{section}
\usepackage{indentfirst}

\begin{document}

\newtheorem{thm}{Theorem}[section]
\newtheorem{cor}[thm]{Corollary}
\newtheorem{prop}[thm]{Proposition}
\newtheorem{conj}[thm]{Conjecture}
\newtheorem{lem}[thm]{Lemma}
\newtheorem{Def}[thm]{Definition}
\newtheorem{rem}[thm]{Remark}
\newtheorem{prob}[thm]{Problem}
\newtheorem{ex}{Example}[section]

\newcommand{\be}{\begin{equation}}
\newcommand{\ee}{\end{equation}}
\newcommand{\ben}{\begin{enumerate}}
\newcommand{\een}{\end{enumerate}}
\newcommand{\beq}{\begin{eqnarray}}
\newcommand{\eeq}{\end{eqnarray}}
\newcommand{\beqn}{\begin{eqnarray*}}
\newcommand{\eeqn}{\end{eqnarray*}}
\newcommand{\bei}{\begin{itemize}}
\newcommand{\eei}{\end{itemize}}

\newcommand{\pa}{{\partial}}
\newcommand{\V}{{\rm V}}
\newcommand{\R}{{\bf R}}
\newcommand{\K}{{\rm K}}
\newcommand{\e}{{\epsilon}}
\newcommand{\tomega}{\tilde{\omega}}
\newcommand{\tOmega}{\tilde{Omega}}
\newcommand{\tR}{\tilde{R}}
\newcommand{\tB}{\tilde{B}}
\newcommand{\tGamma}{\tilde{\Gamma}}
\newcommand{\fa}{f_{\alpha}}
\newcommand{\fb}{f_{\beta}}
\newcommand{\faa}{f_{\alpha\alpha}}
\newcommand{\faaa}{f_{\alpha\alpha\alpha}}
\newcommand{\fab}{f_{\alpha\beta}}
\newcommand{\fabb}{f_{\alpha\beta\beta}}
\newcommand{\fbb}{f_{\beta\beta}}
\newcommand{\fbbb}{f_{\beta\beta\beta}}
\newcommand{\faab}{f_{\alpha\alpha\beta}}

\newcommand{\pxi}{ {\pa \over \pa x^i}}
\newcommand{\pxj}{ {\pa \over \pa x^j}}
\newcommand{\pxk}{ {\pa \over \pa x^k}}
\newcommand{\pyi}{ {\pa \over \pa y^i}}
\newcommand{\pyj}{ {\pa \over \pa y^j}}
\newcommand{\pyk}{ {\pa \over \pa y^k}}
\newcommand{\dxi}{{\delta \over \delta x^i}}
\newcommand{\dxj}{{\delta \over \delta x^j}}
\newcommand{\dxk}{{\delta \over \delta x^k}}

\newcommand{\px}{{\pa \over \pa x}}
\newcommand{\py}{{\pa \over \pa y}}
\newcommand{\pt}{{\pa \over \pa t}}
\newcommand{\ps}{{\pa \over \pa s}}
\newcommand{\pvi}{{\pa \over \pa v^i}}
\newcommand{\ty}{\tilde{y}}
\newcommand{\bGamma}{\bar{\Gamma}}

\font\BBb=msbm10 at 12pt
\newcommand{\Bbb}[1]{\mbox{\BBb #1}}

\newcommand{\qed}{\hspace*{\fill}Q.E.D.}  

\title{The Randers metrics of weakly isotropic scalar curvature}
\author{Xinyue Cheng\footnote{supported by the National Natural Science Foundation of China (No.11871126) and Chongqing Normal University Science Research Fund (No. 17XLB022)}, Yannian Gong}

\maketitle

\begin{abstract}
In this paper, we study the Randers metrics of weakly isotropic scalar curvature.  We  prove that a Randers metric of weakly isotropic scalar curvature must be of isotropic $S$-curvature. Further, we prove that a conformally flat Randers metric of weakly isotropic scalar curvature is either Minkowskian or Riemannian.\\
{\bf Keywords:} Finsler geometry, Randers metric, Ricci curvature tensor,  scalar curvature, $S$-curvature.
\end{abstract}

\section{Introduction}

Randers metrics form  a special and important class of metrics in Finsler geometry. A Randers metric on a manifold $M$ is a Finsler metric in the following form:
\[
F = \alpha +\beta,
\]
where $\alpha =\sqrt{a_{ij}(x)y^iy^j}$ is a Riemannian metric and $\beta = b_i(x) y^i$ is a $1$-form  satisfying $\|\beta_x\|_{\alpha} < 1$ on  $M$. Randers metrics
were first introduced  by physicist G. Randers  in 1941 from the standpoint of general relativity, here the Riemannian metric $\alpha$ denotes the gravitation field and $\beta$ denotes a electromagnetic field. Later on, these metrics were applied to the theory of the electron microscope by R. S. Ingarden in 1957, who first named them Randers metrics. An interesting fundamental fact about Randers metrics is as follows: any Randers metric can be expressed as the solution of Zermelo navigation problem with navigation data $(h, W)$, where $h$ is a Riemannian metric and $W$ is a vector field with $h(x, -W)<1$ on $M$ (\cite{CS}).

Finsler geometry is just the Riemannian geometry without quadratic restriction. Ricci curvature in Finsler geometry is the natural extension of that in Riemannian geometry. However,
there is no unified definition of Ricci curvature tensor in Finsler geometry. Hence, we can find  several different versions of the definition of scalar curvature in Finsler geometry. Here, we adopt the definitions introduced by H. Akbar-Zadeh for Ricci curvature tensor and scalar curvature (\cite{AZ}). For a Finsler metric $F$ on an $n$-dimensional manifold $M$, let $\bf Ric$ be the Ricci curvature of $F$. Then the scalar curvature of $F$ is defined as folows
\be
{\bf r}:=g^{ij}{\bf Ric}_{ij},\label{eqb0}
\ee
where
\beqn
{\bf Ric}_{ij}:=\frac {1}{2}{\bf Ric}_{y^{i}y^{j}}
\eeqn
denote the Ricci curvature tensor and $(g^{ij}):=(g_{ij})^{-1}, \ g_{ij}:=\frac {1}{2}[F^{2}]_{y^{i}y^{j}}$. We say that $F$ is of weakly isotropic scalar curvature if there exists a 1-form $\theta:=\theta_{i}(x)y^{i}$ and a scalar function $\mu (x)$ on $M$ such that
\be
{\bf r}=n(n-1)\left[\frac{\theta}{F}+\mu(x)\right]. \label{wisc}
\ee
In particular, when $\theta=0$, that is, ${\bf r}=n(n-1)\mu(x)$, we say that $F$ is of isotropic scalar curvature. We can find many Finsler metrics of weakly isotropic scalar curvature which are not of isotropic scalar curvature (see Example \ref{ex1} below).

The $S$-curvature ${\bf S} = {\bf S}(x, y)$ is an important non-Riemannian quantity in Finsler geometry which was first introduced by Z. Shen when he studied
volume comparison in Riemann-Finsler geometry (\cite{Sh}). Shen proved that the Bishop-Gromov volume comparison holds for Finsler manifolds with vanishing $S$-curvature. The recent studies show that $S$-curvature plays a very important role in Finsler geometry. In 2014, the first author and M. Yuan verified that a Randers metric of isotropic scalar curvature must be of isotropic $S$-curvature (see {\cite{CY}}). In this paper, we mainly study the Randers metrics of weakly isotropic scalar curvature. Firstly, we obtain the following theorem which generalizes the related result in \cite{CY} mentioned above.

\begin{thm} \label{SCS} \ Let $F = \alpha +\beta$ be a Randers metric on an $n$-dimensional manifold $M$. If $F$ is of weakly isotropic scalar curvature,  ${\bf r}=n(n-1)\left[\frac{\theta}{F}+\mu(x)\right]$,  then $F$ is of  isotropic $S$-curvature.
\end{thm}

The following is an example about Randers metrics of weakly isotropic scalar curvature which arises from \cite{CS1}.
\begin{ex}{\rm (\cite{CS1})}\label{ex1} Let us consider the following Randers metric
\beqn
F&=& \frac{\sqrt{\left(1-|a|^{2}|x|^{4}\right)|y|^{2}+\left(|x|^{2}\langle a, y\rangle- 2\langle a, x\rangle\langle x, y\rangle\right)^{2}}}{1-|a|^{2}|x|^{4}} \\
&& -\frac{|x|^{2}\langle a, y\rangle- 2\langle a, x\rangle\langle x, y\rangle}{1-|a|^{2}|x|^{4}},
\eeqn
where $a$ is a constant vector in ${\bf R}^{n}$ and $\langle ~, \rangle$ denotes the standard inner product in ${\bf R}^{n}$. By direct computation, one can easily verify that $F$ is of weakly isotropic scalar curvature. Precisely, we have
\be
{\bf r}=n(n-1)\left(\frac{\theta}{F}+\mu(x)\right), \ \ \ \theta =\frac{3(n+1)c_{m}y^{m}}{2n}, \label{ex1}
\ee
where $c=\langle a, x\rangle$ and $\mu = 3\langle a, x\rangle^{2}-2|a|^{2}|x|^{2}$, \ $c_{m}:=c_{x^{m}}$. Further, we can also prove that $F$ is of isotropic $S$-curvature,
\[
{\bf S}=(n+1) c F.
\]
Actually, in this case, $F$ is of weakly isotropic flag curvature,
\[
\mathbf{K}=\frac{3 c_{m} y^{m}}{F}+ \mu (x).
\]
Then we can get (\ref{ex1}) by Lemma \ref{wEisc}.
\end{ex}

The study on conformal geometry has a long and venerable history. From the beginning, conformal geometry has played an important role in differential geometry and physical theories. The Weyl theorem shows that the conformal and projective properties of a Finsler space determine the properties of metric completely (see{\cite{Kn}}, {\cite{Ru}} and \cite{BC}).  Undoubtedly, Finsler conformal geometry is an important part of Finsler geometry.  We say two Finsler metrics $F$ and $\bar{F}$ are conformally related if there is a scalar function $\sigma (x)$ on the manifold such that $F=e^{\sigma(x)}\bar{F}$. Further, if $\bar{F}$ is a Minkowskian, the Finsler metric $F$ is called the conformally flat Finsler metric. It is an important topic in Finsler geometry to reveal and characterize deeply geometric structures and properties of conformally flat Finsler metrics.
G. Chen and the first author have proved that a conformally flat weak Einstein $(\alpha, \beta)$-metric must be either a locally Minkowski metric or a Riemannian metric on a
manifold M with the dimension $n \geq 3$ ({\cite {CC}}).  Chen-He-Shen proved that a conformally flat $(\alpha, \beta)$-metric with constant flag curvature is either a locally Minkowskian or Riemannian metric (\cite{GQZ}). On the other hand, Cheng-Yuan proved that a conformally flat non-Riemannian Randers metric of isotropic scalar curvature must be locally Minkowskian when the dimension $n \geq 3$ (\cite{CY}). Further,
B. Chen and K. Xia studied a class of conformally flat polynomial $(\alpha, \beta)$-metrics in the form $F=\alpha\left(1+\Sigma^m_{j=2}a_{j}(\frac {\beta}{\alpha})^{j}\right)$ with $m\geq 2$.  They proved that,
if such a conformally flat $(\alpha, \beta)$-metric F is of weaklly isotropic scalar curvature, then it  must have zero salar curvature. Moreover, if $a_{m-1}a_{m}\neq 0$, then $F$ must be either locally Minkowkian or Riemannian when the dimension  $n\geq 3$ (\cite{BK}). When $m=1$, that is, when $F$ is a Randers metric, Chen-Xia have not confirmed that the same conclusion still holds. Therefore,  it is a natural problem  to characterize conformally flat Randers metrics of weakly isotropic scalar curvature.   We have got the following theorem.

\begin{thm} \label{CWS}\ Let $F = \alpha +\beta$ be a conformally flat non-Riemannian Randers metric on an n-dimensional manifold $M$ with $n\geq 2$. If $F$ is of weakly isotropic scalar curvature, that is, ${\bf r}=n(n-1)[\frac{\theta}{F}+\mu(x)]$, then $F$ must be locally Minkowkian.
\end{thm}

\section{Preliminaries}

 Let $M$ be an $n$-dimensional smooth manifold and $(x^{i}, y^{i})$ denote the local coordinates of point $(x,y)$ on the tangent boudle $TM$ with $y=y^{i}\frac{\pa}{\pa x^{i}}\in T_{x}M$. Let $F$ be a Finsler metric on $M$ and $g_{y}=g_{kl}(x, y)dx^{k}\otimes dx^{l}$ be the fundamental tensor of $F$, where $g_{kl}:=\frac {1}{2}[F]^{2}_{y^{k}y^{l}}$. The  geodesic coefficients of $F$ are given by
 \be
 G^{k}=\frac {1}{4}g^{kl}\{[F]^{2}_{x^{m}y^{l}}y^{m}-[F^{2}]_{x^{l}}\},
 \ee
where $(g^{kl}):=(g_{kl})^{-1}$.  For any $x\in M$ and $y\in T_xM\backslash \left\{0\right\}$, the  Riemann curvature ${\bf R}_{y}:=R^{i}_{\ k}(x, y)\frac {\partial}{\partial x^{i}}$ $\otimes dx^{k}$ is defined by
 \be
 R^{i}_{\ k}(x, y):=2G^{i}_{x^{k}}-G^{i}_{x^jy^k}y^{j}+2G^{j}G^{i}_{y^{j}y^{k}}-G^{i}_{y^{j}}G^{j}_{y^{k}}.
 \ee
The  Ricci curvature of Finsler metric $F$ is defined as the trace of Riemann curvature , that is
 \be
{\bf Ric}(x, y):=R^{m}_{\ m}(x, y).
 \ee
 It is not difficult to see that Ricci curvature is a positive homogeneous function of degree two in y. Further, the  Ricci curvature tensor is given by
 \be
 {\bf Ric}_{ij}:=\frac {1}{2}{\bf Ric}_{y^{i}y^{j}}.\label{Ricci def}
 \ee
 One can get ${\bf Ric}(x, y)={\bf Ric}_{ij}y^{i}y^{j}$ by the homogeneity of ${\bf Ric}$. A Finsler metric $F$ is called a  weak Einstien metric, if there exists a 1-form $\xi= \xi_{i}(x)y^{i}$ and a scalar function $\mu=\mu(x)$ on $M$ such that
 \be
 {\bf Ric}=(n-1)\left(\frac {3 \xi}{F}+\mu\right)F^{2} . \label{wE}
 \ee
 In particular, if $\xi =0$, that is, ${\bf Ric}=(n-1)\mu F^{2}$, $F$ is called an  Einstein metric.

 The scalar curvature of a Finsler metric $F$ introduced by Akbar-Zadeh is defined by (1.1), that is, ${\bf r}:=g^{ij}{\bf Ric}_{ij}$.  The following lemma is natural and important.
 \begin{lem}\label{wEisc}
 Assume that $F$ is a weak Einstein Finsler metric satisfying (\ref{wE}). Then $F$ must be of weakly isotropic scalar curvature satisfying
 \be
 {\bf r}= n(n-1)\left(\frac{\theta}{F}+\mu\right), \ \ \theta = \frac{3(n+1)}{2n}\xi. \label{wEWIS}
 \ee
 \end{lem}

The distortion $\tau$ of a Finsler metric $F$ is defined by
\[
\tau(x, y):=\ln\frac {\sqrt{det(g_{ij}(x, y))}}{\sigma_{BH}(x)},
\]
 where
 \[
 \sigma_{BH}:=\frac{Vol({\bf B}^{n}(1))}{Vol\{(y^{i})\in R^{n}|F(x, y^{i}\frac{\partial}{\partial x^{i}})<1\}}
  \]
 is Busemann-Hausdorff volume coefficient. $\tau=0$ if and only if  Finsler metric $F$ is Riemannian (\cite{SSZ}). $S$-curvature ${\bf S}$ of $F$ characterizes the change rate of distortion  $\tau$ along geodesics, that is,
 \[
 {\bf S}(x,y):=\tau_{|m}(x, y)y^{m}.
 \]
In local coordinate system, the $S$-curvature of $F$ can be expressed as
 \[
 {\bf S}(x, y)=\frac{\pa G^{m}}{\pa y^{m}}-y^{m}\frac{\pa}{\pa x^{m}}\left[\ln\sigma_{BH}(x)\right].
 \]
 We say that a Finsler metric $F$ is of  isotropic $S$-curvature, if there is a scalar function $c(x)$ on $M$ such that
 \be
 {\bf S}(x, y)=(n-1)c(x)F(x, y).
 \ee

The mean Cartan torsion ${\bf I}_{y}=I_{i}(x, y) d x^{i}: T_{x} M \rightarrow {\bf R}$ is defined by
$$
I_{i}:=g^{j k} C_{i j k},
$$
where $C_{ijk}$ denote the Cartan torsion of $F$.  It is easy to check that  $ I_{i}=\tau_{y^{i}}$.

\vskip 2mm

For a Randers metric $F=\alpha+\beta$,  we have the following
\be\label{gij}
g^{ij}=\frac {\alpha}{F}a^{ij}-\frac {\alpha}{F^{2}}(b^{i}y^{j}+b^{j}y^{i})+\frac {b^{2}\alpha+\beta}{F^{3}}y^{i}y^{j},
\ee
where $b:=\|\beta\|_{\alpha}$ denotes the norm of $\beta$ with respect to $\alpha$ (\cite{SSZ}). Let
\[
r_{ij}:=\frac {1}{2}(b_{i;j}+b_{j;i}),  \ \ \ \  s_{ij}:=\frac {1}{2}(b_{i;j}-b_{j;i}),\ \ \ \ \ \ \ \ \ \\
\]
\[
e_{ij}:=r_{ij}+s_{i}b_{j}+s_{j}b_{i},\ \ \ \ \ \ \ \ \ \ \ \ \ \ \ \ \ \ \ \ \ \\
\]
\[
w_{ij}:=r_{im}r^{m}\!_{j},\ \ t_{ij}:=s_{im}s^{m}\!_{j},\ \ q_{ij}:=r_{im}s^{m}\!_{j}, \ \ \ \ \ \ \  \\
\]
\[
r^{i}_{\ j}:=a^{im}r_{mj},\ s^{i}_{\ j}:=a^{im}s_{mj}, \ t^{i}_{\ j}:=a^{im}t_{mj},\ q^{i}_{\ j}:=a^{im}q_{mj},\\
\]
\[
r_{i}:=b^{m}r_{mi},\ s_{i}:=b^{m}s_{mi},\ t_{i}:=b^{m}t_{mi},\ q_{i}:=b^{m}q_{mi},\ \ \\
\]
\[
r:=b^{i}b^{j}r_{ij},\ \ \ t:=b^{i}t_{i},\ \ \ p_{i}:=r_{im}s^{m},
\]
where `` ; " denotes the covariant derivative with respect to $\alpha$. Besides, put $ r_{00}:=r_{ij}y^{i}y^{j},\ e_{00}:=e_{ij}y^{i}y^{j},\ q_{00}:=q_{ij}y^{i}y^{j},\ s_{0}:=s_{i}y^{i}$,  etc..

Further, the Ricci curvature of Randers metric $F=\alpha +\beta$ is given by
\be
{\bf Ric}={}^{\alpha}{\bf Ric}+(2\alpha s^{m}\!_{0;m}-2t_{00}-\alpha^{2}t^{m}_{\ m})+(n-1)\Xi,\label{Ricci curvature}
\ee
where ${}^{\alpha}{\bf Ric}$ denotes the Ricci curvature of $\alpha$ and
\be
\Xi:=\frac {2\alpha}{F}(q_{00}-\alpha t_{0})+\frac {3}{4F^{2}}(r_{00}-2\alpha s_{0})^{2}-\frac
 {1}{2F}(r_{00;0}-2\alpha s_{0;0}). \label{Xi}
\ee
Besides, the mean Cartan tensor ${\bf I}=I_{i}dx^{i}$ of $F=\alpha +\beta$ is given by
\be\label{Ii}
I_{i}=\frac {n+1}{2F}(b_{i}-\frac {\beta y_{i}}{\alpha^{2}}),
\ee
where $y_{i}:=a_{ij}y^{j}$. For related details, see \cite{CS} or \cite{SSZ}.

The following  lemma is very important for the proof of Theorem \ref{SCS}.
\begin{lem}{\rm (\cite{CS})}\label{e00}
Let $F = \alpha +\beta$ be a Randers metric on an $n$-dimensional manifold $M$. Then $F$ is of isotropic $S$-curvature, ${\bf S}=(n+1)c(x)F$ if and only if
\be
e_{00}=2c(x)(\alpha^{2}-\beta^{2}),
\ee
where  $c=c(x)$ is a scalar function on $M$.
\end{lem}

\section{The scalar curvature of Randers metric}\label{section3}

In \cite{CY}, the first author and M. Yuan have obtained the formula of the scalar curvature of Randers metrics.  In this section, we will further improve and optimize this formula.  For readers' convenience, we will derive the formula of the scalar curvature of Randers metrics step by step.

By  (\ref{gij}), (\ref{Ricci curvature}) and the definition of scalar curvature, we can get
\be
{\bf r}:=g^{ij}{\bf Ric}_{ij}= {}^{\alpha}{\bf Ric}_{ij}g^{ij}+\frac {1}{2}E_{ij}g^{ij}+\frac {1}{2}(n-1)\Xi_{ij}g^{ij},\label{scRanders}
\ee
where  ${}^{\alpha}{\bf Ric}_{ij}$ denote Ricci curvature tensor of $\alpha$ and
\beqn
&& E:= 2\alpha s^{m}_{\ 0;m}-2t_{00}-\alpha^{2}t^{m}_{\ m},\\
&& E_{ij}:= E_{y^{i}y^{j}},\ \ \ \ \Xi_{ij}:=\Xi_{y^{i}y^{j}}.
\eeqn
Next, we will compute each term on the right side of (\ref {scRanders}). Firstly, we get the following
\be\label{term 1}
{}^{\alpha}{\bf Ric}_{ij}g^{ij}=\frac {\alpha}{F}{\bf r}_{\alpha}-\frac{2\alpha}{F^{2}}\ {}^{\alpha}{\bf Ric}_{ij}b^{i}y^{j}+\frac {b^{2}\alpha+\beta}{F^{3}}\ {}^{\alpha}{\bf Ric},
\ee
where ${\bf r}_{\alpha}$ denotes the scalar curvature of $\alpha$. Moreover, we obtain the following
\beq
E_{ij}g^{ij}&=& \frac{2\alpha}{F}\left[\frac{n+1}{\alpha}s^{m}_{\ 0;m}-(n+2)t^{m}_{\ m}\right] \nonumber\\
 &&-\frac{4\alpha}{F^{2}}\big[s^{m}_{\ 0;m}s+\alpha b^{i} s^{m}_{\ i;m}-2t_{0}-\beta t^{m}_{\ m}\big]+2\frac {b^{2}\alpha+\beta}{F^{3}}E, \label{Eij}
\eeq
where $s:=\beta / \alpha$.

In order to obtian $\Xi_{ij}g^{ij}$, we rewrite (\ref{Xi}) as the following
\be
 \Xi:=\Xi^{1}+\Xi^{2}+\Xi^{3}, \label{reXi}
\ee
where,
\[
\Xi^{1}:=2\left(\frac {D}{F}\right), \ \ \Xi^{2}:=\frac {3}{4}\left(\frac {A^{2}}{F^{2}}\right), \ \ \Xi^{3}:= -\frac {1}{2}\left(\frac {B}{F}\right),
\]
and
\[
A:=r_{00}-2\alpha s_{0},\ \ \ B:=r_{00;0}-2\alpha s_{0;0},
\]
\[
 D:=\alpha D_{1}=\alpha(q_{00}-\alpha t_{0}), \ \ \ D_{1}:= q_{00}-\alpha t_{0}.
\]

From  (\ref{reXi}), one can write
\be
 \Xi_{ij}:=\Xi^{1}_{ij}+\Xi^{2}_{ij}+\Xi^{3}_{ij},\label{XXi}
\ee
where
\beqn
\Xi^{1}_{ij}&:=& \Big[2\frac{\alpha}{F}(q_{00}-\alpha t_{0})\Big]_{y^{i}y^{j}}=2\Big(\frac{D}{F}\Big)_{y^{i}y^{j}},\\
\Xi^{2}_{ij}&:=& \Big[\frac{3}{4F^{2}}(r_{00}-2\alpha s_{0})^{2}\Big]_{y^{i}y^{j}}=\Big(\frac{3}{4}\Big)\Big(\frac{A^{2}}{F^{2}}\Big)_{y^{i}y^{j}},\\
\Xi^{3}_{ij}&:=& \Big[-\frac{1}{2F}(r_{00;0}-2\alpha s_{0;0})\Big]_{y^{i}y^{j}}=-\Big(\frac{1}{2}\Big)\Big(\frac{B}{F}\Big)_{y^{i}y^{j}}.
\eeqn
By a series direct computations, we can get the following
\beq
\Xi^{1}_{ij}g^{ij}&=& \frac{2}{F}\left \{\frac{\alpha}{F}\Big[(n+3)\frac{D_{1}}{\alpha}+2\alpha q^{m}\!_{m}-(n+1)t_{0}\Big]\right. \nonumber \\
&& \left.-\frac{4\alpha}{F^{2}}\big[sD_{1}+\alpha (q_{00\cdot i}b^{i}-st_{0}-\alpha t)\big]+6D\frac{b^{2}\alpha+\beta}{F^{3}}\right\} \nonumber\\
&& -\frac{4}{F^{2}}\left\{\frac{\alpha}{F}\big[(3+s)D_{1}+\alpha (q_{00\cdot i}b^{i}-st_{0}-\alpha t)\big]\right. \nonumber\\
&& \left. -\frac{\alpha}{F^{2}}\big[ F(sD_{1}+\alpha (q_{00\cdot i}b^{i}-st_{0}-\alpha t))+3D(b^{2}+s)\big]+ 3D\frac{b^{2}\alpha+\beta}{F^{2}}\right\} \nonumber\\
&& -2(n-1)\frac{D}{F^{3}}+4\frac{D}{F^{3}}\Big[\frac{\alpha}{F}(1-b^{2})+\frac{b^{2}\alpha+\beta}{F}\Big] \nonumber\\
& =& \frac{1}{2F^{5}}\left\{4F^{3}\big[(n+3)D_{1}-(n+1)\alpha t_{0}+2\alpha^{2}q^{m}_{\ m}\big]\right. \nonumber\\
&& -16F^{2}\big[\beta D_{1}+\alpha^{2}q_{00\cdot i}b^{i}-\alpha\beta t_{0}-\alpha^{3}t\big]- 4F^{2}(n+3)D \nonumber\\
&& \left.+24F(b^{2}\alpha+\beta)D\right\}, \label{Xieq1}
\eeq
\beq
\Xi^{2}_{ij}g^{ij}&=&\frac{6}{F^{2}}\left\{\frac{\alpha}{F}\Big[w_{00}-2\frac{r_{00}s_{0}}{\alpha}-2\alpha p_{0}+3s^{2}_{0}-\alpha^{2}t\Big]\right.\nonumber\\
&& \left.-\frac{2\alpha}{F^{2}}A(r_{0}-ss_{0})+\frac{b^{2}\alpha+\beta}{F^{3}}A^{2}\right\}+\frac{3A}{F^{2}}\left\{\frac{\alpha}{F}\big[r^{m}_{\ m}-(n+1)\frac{s_{0}}{\alpha}\big]\right. \nonumber\\
&& \left.-\frac{2\alpha}{F^{2}}(r_{0}-ss_{0})+\frac{b^{2}\alpha+\beta}{F^{3}}A\right\}-\frac{12A}{F^{3}}\left\{\frac{\alpha}{F}\Big[\frac{r_{00}}{\alpha}+r_{0}-(2+s)s_{0}\Big] \right. \nonumber\\
&& \left.-\frac{\alpha}{F^{2}}\Big[F(r_{0}-ss_{0})+A(s+b^{2})\Big]+\frac{b^{2}\alpha+\beta}{F^{2}}A\right\}-\frac{3(n-1)A^{2}}{2F^{4}} \nonumber\\
&& +\frac{9A^{2}}{2F^{4}}\Big[\frac{\alpha}{F}(1-b^{2})+\frac{b^{2}\alpha+\beta}{F}\Big] \nonumber\\
&=& \frac{1}{2F^{5}}\left\{12F^{2}(\alpha w_{00}-2r_{00}s_{0}-2\alpha^{2}p_{0}+3\alpha s^{2}_{0}-\alpha^{3}t)\right.\nonumber\\
&& +6F^{2}A\big[\alpha r^{m}\!_{m}-(n+1)s_{0}\big]-12FA(\alpha r_{0}-\beta s_{0}) \nonumber\\
&& -24FA\big[r_{00}+\alpha r_{0}-(2\alpha+\beta)s_{0}\big]-3(n-4)FA^{2} \nonumber\\
&& \left.+18A^{2}(b^{2}\alpha+\beta)\right\} \label{Xieq2}
\eeq
and
\beq
\Xi^{3}_{ij}g^{ij}&=&\frac{1}{2F}\left\{-\frac{2\alpha}{F}\big[2r^{m}_{\ 0;m}+r^{m}_{\ m;0}-(n+3)\frac{s_{0;0}}{\alpha}-2\alpha s^{m}_{\ ;m}\big]\right. \nonumber\\
&& \left.+\frac{4\alpha}{F^{2}}\big[r_{00;0\cdot i}b^{i}-2ss_{0;0}-2\alpha s_{0;0\cdot i}b^{i}\big]-2\frac{b^{2}\alpha+\beta}{F^{3}}B\right\} \nonumber\\
&& +\frac{1}{F^{2}}\left\{\frac{\alpha}{F}\big[\frac{3r_{00;0}}{\alpha}-2(3+s)s_{0;0}+r_{00;0\cdot i}b^{i}-2\alpha s_{0;0\cdot i}b^{i}\big]\right. \nonumber\\
&& \left.-\frac{\alpha}{F^{2}}\big[F(r_{00;0\cdot i}b^{i}-2ss_{0;0}-2\alpha s_{0;0\cdot i}b^{i})+3B(s+b^{2})\big]\right\} \nonumber\\
&& +\frac{n-1}{2F^{3}}B-\frac{ B}{F^{4}}\alpha(1-b^{2})\nonumber\\
&=&\frac{1}{2F^{5}}\left\{-2F^{3}\big[\alpha r^{m}_{\ m;0}+2\alpha r^{m}_{\ 0;m}-2\alpha^{2}s^{m}_{\ ;m}-(3+n)s_{0;0}\big]\right.\nonumber\\
&& +F^{2}\big[4\alpha r_{00;0\cdot i}b^{i}-8\alpha^{2}s_{0;0\cdot i}b^{i}-8\beta s_{0;0}+(n+3)B \big] \nonumber\\
&& \left.-6F(b^{2}\alpha+\beta)B\right\}, \label{Xieq3}
\eeq
where `` $\cdot i$ " denotes the partial derivative with respect to $y^{i}$.

Now, plugging (\ref{term 1}), (\ref{Eij}) and (\ref{Xieq1})-(\ref{Xieq3}) into (\ref{scRanders}), one obtains the formula of the scalar curvature for Randers metric $F=\alpha +\beta$ as follows
\be
{\bf r}=\frac {\alpha}{F}{\bf r}_{\alpha}+\frac {1}{4F^{5}}\Big\{\Sigma_{1}+\alpha\Sigma_{2}\Big\},  \label{r1}
\ee
where $\Sigma_{1}$ and $\Sigma_{2}$ are both polynomials in $y$.  Concretely, we have the following expressions:
\[
 \begin{aligned}
\Sigma_{1}:=&\left\{\big[-4(2b^{2}+4n+7)t^{m}_{\ m}-24b^{i}s^{m}_{\ i;m}+4(n-1)(6q^{m}_{\ m}+3s^{m}_{\ ;m}+2t)\big]\beta\right.\\
&-8\ {}^{\alpha}{\bf Ric}_{ij}b^{i}y^{j}+4(2b^{2}+n+1)s^{m}_{\ 0;m}-4(6b^{2}n-6b^{2}+n^{2}-5)t_{0}\\
&\left.-2(n-1)(6r^{m}_{\ m}s_{0}+8q_{00\cdot i}b^{i}+r^{m}_{\ m;0}+2r^{m}_{\ 0;m}+4s_{0;0\cdot i}b^{i}+12p_{0})\right\}\alpha^{4}\\
&+\left\{-4\big[(3+4n)t^{m}_{\ m}+2b^{i}s^{m}_{\ i;m}-(n-1)(2q^{m}_{\ m}+s^{m}_{\ ;m})\big]\beta^{3}\right.\\
&+\big[-24 {}^{\alpha}{\bf Ric}_{ij}b^{i}y^{j}+8(b^{2}+3n+2)s^{m}_{\ 0;m}-4(n+1)(5n-11)t_{0}\\
&-2(n-1)(6r^{m}_{\ m}s_{0}+8q_{00\cdot i}b^{i}+3r^{m}_{\ m;0}+6r^{m}_{\ 0;m}+4s_{0;0\cdot i}b^{i}+12p_{0})\big]\beta^{2}\\
&+\big[4(2b^{2}+1)~ {}^{\alpha}{\bf Ric}-8(2b^{2}+1)t_{00}+4(n-1)(6b^{2}+n+5)q_{00}\\
&+12(n-1)(3n-4)s_{0}^{2}+2(n-1)(6b^{2}+n+5)s_{0;0}+4(n-1)(3r^{m}_{\ m}r_{00}\\
&+18r_{0}s_{0}+2r_{00;0\cdot i}b^{i}+6w_{00})\big]\beta-6(n-1)\big[(12b^{2}+3n-19)s_{0}+6r_{0}\big] r_{00} \\
& \left.-(n-1)(6b^{2}-n-3)r_{00;0} \right\}\alpha^{2}+4(n-1)s^{m}_{\ 0;m}\beta^{4}+\big[4 {}^{\alpha}{\bf Ric}-8t_{00} \\
&+2(n-1)^{2}(2q_{00}+s_{0;0})\big]\beta^{3}+(n-1)\big[-6(n-1)s_{0}r_{00}+(n-3)r_{00;0}\big]\beta^{2}\\
&+3(n-1)(n-6)r_{00}^{2}\beta
 \end{aligned}
\]
and
\[
 \begin{aligned}
\Sigma_{2}:=&\big[-4(b^{2}+n+2)t^{m}_{\ m}-8b^{i}s^{m}_{\ i;m}+4(n-1)(2q^{m}_{\ m}+s^{m}_{\ ;m}+t)\big]\alpha^{4}\\
&+\left\{[-4(b^{2}+6n+8)t^{m}_{\ m}-24b^{i}s^{m}_{\ i;m}+4(n-1)(6q^{m}_{\ m}+3s^{m}_{\ ;m}+t)]\beta^{2}\right.\\
&+[-24 {}^{\alpha}{\bf Ric}_{ij}b^{i}y^{j}+16(b^{2}+n+1)s^{m}_{\ 0;m}-2(n-1)(12r^{m}_{\ m}s_{0}+24p_{0}\\
&+16q_{00\cdot i}b^{i}+3r^{m}_{\ m;0}+6r^{m}_{\ 0;m}+8s_{0;0\cdot i}b^{i})-8(3b^{2}n-3b^{2}+2n^{2}-8)t_{0}\big]\beta\\
&+4\ {}^{\alpha}{\bf Ric}b^{2}-8b^{2}t_{00}+2(n-1)\big[12(3b^{2}+n-4)s_{0}^{2}+36s_{0}r_{0}+6b^{2}s_{0;0}\\
&\left.+12b^{2}q_{00}+3r^{m}_{\ m}r_{00}+2r_{00;0\cdot i}b^{i}+6w_{00}]\right\}\alpha^{2}-4nt^{m}_{\ m}\beta^{4}+2\big[8ns^{m}_{\ 0;m}\\
&+(n-1)(2r^{m}_{\ 0;m}+r^{m}_{\ m;0})-4 {}^{\alpha}{\bf Ric}_{ij}b^{i}y^{j}-4n(n-3)t_{0}\big]\beta^{3}\\
&+\left\{4(b^{2}+2)({}^{\alpha}{\bf Ric}-2t_{00})+2(n-1)\big[6(n-2)s_{0}^{2}+2(n+2)s_{0;0}+6w_{00} \right. \\
&\left.+3r^{m}_{\ m}r_{00}+4(n+4)q_{00}+2r_{00;0\cdot i}b^{i}\big]\right\}\beta^{2}-2(n-1)\left\{3b^{2}r_{00;0}-nr_{00;0}\right.\\
&\left.+6[2(n-2)s_{0}-3r_{0}]r_{00}\right\}\beta+3(n-1)(6b^2+n-12)r_{00}^{2}.
 \end{aligned}
\]
It is obvious that $\Sigma_{1}$ and  $\Sigma_{2}$ are homogeneous polynomials of degree 5 and 4 in $y$, respectively.

\begin{rem} We have totally corrected some errors occurred in the formulas of  $\Sigma_{1}$ and  $\Sigma_{2}$  in \cite{CY}. We must point out that, in the proofs of main theorems in \cite{CY}, the authors just used the facts that $\Sigma_{1}$ and  $\Sigma_{2}$ are homogeneous polynomials of dgree 5 and 4 in $y$, respectively. Hence, the main results in \cite{CY} are still true.
\end{rem}

\section{Proof of Theorems}

In this section, we are going to prove Theorem \ref{SCS} and Theorem \ref{CWS}.
\vskip 2mm

{\bf Proof of Theorem \ref{SCS}.} \ Let $F=\alpha+\beta$ be a Randers metric on an $n$-dimentional manifold $M$ with weakly isotropic scalar curvature.

Firstly, note that
\[
 e_{ij}:=r_{ij}+b_{i}s_{j}+b_{j}s_{i}.
\]
We have
\be
r_{00}=e_{00}-2\beta s_{0} \label{r00}
\ee
and
\be
r_{00;0}=e_{00;0}-2(\beta s_{0;0}+s_{0}e_{00}-2\beta s_{0}^{2}). \label{r00;0}
\ee
Then plugging (\ref{r00}), (\ref{r00;0}) into (\ref{r1}) and multiplying  ($\ref{r1}$) by $4F^{5}$, one gets
\be
4F^{5}{\bf r}=\Gamma_{1}+\alpha\Gamma_{2},\label{r2}
\ee
where $\Gamma_{1}$ and $\Gamma_{2}$ are homogeneous polynomials of degree 5 and 4 in $y$, respectively, which have the following expressions:

\beq
\Gamma_{1}&:=&\left\{\big[16{\bf r}_{\alpha}-4(2b^{2}+4n+7)t^{m}_{\ m}-24b^{i}s^{m}_{\ i;m}+4(n-1)(6q^{m}_{\ m}+3s^{m}_{\ m}+2t)\big]\beta\right.\nonumber\\
&& -8 {}^{\alpha}{\bf Ric}_{ij}b^{i}y^{j}+4(2b^{2}+n+1)s^{m}_{\ 0;m}-2(n-1)\big[6r^{m}_{\ m}s_{0}+8q_{00\cdot i}b^{i} \nonumber\\
&& \left. +r^{m}_{\ m;0}+2r^{m}_{\ 0;m}+4s_{0;0\cdot i}b^{i}+12p_{0}\big]-4(6b^{2}n-6b^{2}+n^{2}-5)t_{0}\right\}\alpha^4 \nonumber\\
&& + \left\{\big[16{\bf r}_{\alpha}-4(4n+3)t^{m}_{\ m}-8b^{i}s^{m}_{\ i;m}+4(n-1)(2q^{m}_{\ m}+s^{m}_{\ m})\big]\beta^3\right. \nonumber\\
&& +\big[-2(n-1)(18r^{m}_{\ m}s_{0}+ 8q_{00\cdot i}b^{i}+3r^{m}_{\ m;0}+6r^{m}_{\ 0;m}+4s_{0;0\cdot i}b^{i}+12 p_{0}) \nonumber\\
&& -24 {}^{\alpha}{\bf Ric}_{ij}b^{i}y^{j}+8(b^{2}+3n+2)s^{m}_{\ 0;m}-4(n+1)(5n-11)t_{0}\big]\beta^2 \nonumber \\
&& +\Big[4(2b^2+1)({}^{\alpha}{\bf Ric}-2t_{00})+4(n-1)\big[(6b^{2}+n+5)q_{00}+2r_{00;0\cdot i}b^{i} \nonumber\\
&& +(6b^{2}+1)s_{0;0}+3r^{m}_{\ m}e_{00}+36s_{0}r_{0}+(30b^{2}+19n-66)s_{0}^{2}+6w_{00}\big]\Big]\beta \nonumber\\
&& \left. -(n-1)(6b^{2}-n-3)e_{00;0}-4(n-1)\big[9r_{0}+(15b^{2}+5n-27)s_{0}\big]e_{00}\right\}\alpha^{2} \nonumber\\
&& +4(n-1)s^{m}_{\ 0;m}\beta^{4}+\left\{4 {}^{\alpha}{\bf Ric}-8t_{00}+4(n-1)\big[(n-1)q_{00}+s_{0;0}\right. \nonumber\\
&& \left.+(7n-24)s_{0}^{2}\big]\right\}\beta^{3}+(n-1)\big[(n-3)e_{00;0}-4(5n-21)s_{0}e_{00}\big]\beta^{2} \nonumber\\
&& +3(n-1)(n-6)e_{00}^{2}\beta  \label{Gma1}
\eeq
and
\beq
\Gamma_{2}&:=&\left\{4{\bf r}_{\alpha}-4(b^{2}+n+2)t^{m}_{\ m}-8b^{i}s^{m}_{\ i;m}+4(n-1)(2q^{m}_{\ m}+s^{m}_{\ m}+t)\right\}\alpha^{4} \nonumber\\
&& +\left\{\big[24{\bf r}_{\alpha} -4(b^{2}+6n+8)t^{m}_{\ m}-24b^{i}s^{m}_{\ i;m} +4(n-1)(6q^{m}_{\ m}+t +3s^{m}_{\ m})\big]\beta^2 \right. \nonumber\\
&&+\big[-24 {}^{\alpha}{\bf Ric}_{ij}b^{i}y^{j}+16(b^{2}+n+1)s^{m}_{\ 0;m} \nonumber\\
&& -2(n-1)(18r^{m}_{\ m}s_{0}+16q_{00\cdot i}b^{i}+3r^{m}_{\ m;0}+6r^{m}_{\ 0;m}+8s_{0;0\cdot i}b^{i}+24p_{0}) \nonumber\\
&& -8(3b^{2}n-3b^{2}+2n^{2}-8)t_{0}\big]\beta+4 b^{2}\ {}^{\alpha}{\bf Ric}- 8b^{2}t_{00}+2(n-1)\big[12b^{2}q_{00} \nonumber\\
&& \left. +2r_{00;0\cdot i}b^{i}+6b^{2}s_{0;0}+3r^{m}_{\ m}e_{00}+36s_{0}r_{0}+12(3b^{2}+n-4)s_{0}^{2} +6w_{00}\big]\right\}\alpha^{2} \nonumber\\
&& -4(nt^{m}_{\ m}-{\bf r}_{\alpha})\beta^{4}+\big[-2(n-1)(6r^{m}_{\ m}s_{0}+r^{m}_{\ m;0}+2r^{m}_{\ 0;m}) \nonumber\\
&& -8{}^{\alpha}{\bf Ric}_{ij}b^{i}y^{j}+16ns^{m}_{\ 0;m}-8n(n-3)t_{0}\big]\beta^{3}+\left\{4(b^{2}+2)({}^{\alpha}{\bf Ric}-2t_{00})\right. \nonumber\\
&& +2(n-1)\big[4(n+2)q_{00}+2r_{00;0\cdot i}b^{i}+2(3b^{2}+2)s_{0;0}+3r^{m}_{\ m}e_{00}+3s_{0}r_{0} \nonumber\\
&& \left.+4(6b^{2}+10n-33)s_{0}^{2}+6w_{00}\big]\right\}\beta^{2}- 2(n-1)\left\{(3b^{2}-n)e_{00;0}\right. \nonumber\\
&& \left.+\big[18r_{0}+2(15b^{2}+10n-48)s_{0}\big] e_{00}\right\}\beta+3(n-1)(6b^{2}+n-12)e_{00}^{2}.\label{Gma2}
\eeq
By (\ref{Gma1}) and (\ref{Gma2}), we can find that
\be
\Gamma_{2}\beta-\Gamma_{1}=-18(1-b^{2})\beta e_{00}^{2}+(\alpha^{2}-\beta^{2})K_{000},\label{G2-G1}
\ee
where $K_{000}$ is a homogeneous polynomial of degree 3 in $y$.

On the other hand, by the assumption, Randers metric of $F=\alpha +\beta$ is of weakly isotropic scalar curvature, that is,
\be
{\bf r}=n(n-1)\left[\frac{\theta}{F}+\mu(x)\right]. \label{r22}
\ee
Then, multipying  ($\ref{r22}$) by $4F^{5}$, one has
\be \label{r3}
4F^{5}{\bf r}=4n(n-1)(\Pi_{1}+\alpha\Pi_{2}),
\ee
where, $\Pi_{1}$ and $\Pi_{2}$ are homogeneous polynomials of degree 5 and 4 in $y$, respectively, which are expressed as
\beqn
\Pi_{1}&:=&(5\mu\beta+\theta)\alpha^{4}+(10\mu\beta^{3}+6\theta\beta^{2})\alpha^{2}+\mu\beta^{5}+\theta\beta^{4},\\
\Pi_{2}&:=& \mu\alpha^{4}+(10\mu\beta^{2}+4\theta\beta^{2})\alpha^{2}+5\mu\beta^{4}+4\theta\beta^{3}.
\eeqn
Further, we can get
\be
\Pi_{2}\beta-\Pi_{1}=-(\alpha^{2}-\beta^{2})(4\mu\beta\alpha^{2}+4\mu\beta^{3}+\theta\alpha^{2}+3\theta\beta^{2}). \label{r5}
\ee
Comparing (${\ref {r2}}$) and (${\ref {r3}}$), we obtian the following
\be \label{r4}
\Gamma_{1}=4n(n-1)\Pi_{1}, \ \ \Gamma_{2}=4n(n-1)\Pi_{2}.
\ee
Therefore, $\Gamma _{2}\beta - \Gamma _{1}= 4n(n-1)(\Pi_{2}\beta-\Pi_{1})$. From (${\ref {G2-G1}}$) and (${\ref {r5}}$), one obtains
\beqn
&& -18(1-b^{2})\beta e_{00}^{2}+(\alpha^{2}-\beta^{2})K_{000}\\
&&= -4n(n-1)(\alpha^{2}-\beta^{2})(4\mu\beta\alpha^{2}+4\mu\beta^{3}+\theta\alpha^{2}+3\theta\beta^{2}),
\eeqn
which is equivalent to
\be
\-18(1-b^{2})\beta e_{00}^{2}=(\alpha^{2}-\beta^{2})[K_{000}+4n(n-1)(4\mu\beta\alpha^{2}+4\mu\beta^{3}+\theta\alpha^{2}+3\theta\beta^{2})].
\ee
Because  $\alpha^{2}-\beta^{2}$ is an irreducible polynomial in $y$ and $1-b^{2}>0$, we know that $e_{00}$ must be divided by $\alpha^{2}-\beta^{2}$. That is, there exists a scalar function $c(x)$ on $M$ such that
\be
e_{00}=2c(x)(\alpha^{2}-\beta^{2}).
\ee
By Lemma {\ref {e00}}, $F$ is of isotropic $S$-curvature. \qed

\vskip 2mm

In order to prove Theorem \ref{CWS}, we prove the following lemma firstly.

\begin{lem}\label{lmS}
Let $F = \alpha +\beta$ be a non-Riemannian Randers metric with isotropic $S$-curvature on an n-dimensional ($n\geq 2$) manifold $M$. If there is scalar function $\sigma =\sigma (x)$ on $M$ such that $\bar {F}:=e^{\sigma (x)}F$ is of weakly isotropic scalar curvature, then $\sigma$ must be a constant.
\end{lem}

{\bf Proof.} \ By the assumption, $\bar{F}$ is conformally related to $F$, $\bar{F} =e^{\sigma(x)}F$. Then we have the following equality ({\cite {BC}}).
\be
\bar{\bf S} ={\bf S}+F^{2}\sigma^{r}I_{r}, \label{S1}
\ee
where $\sigma^{r}:=g^{rm}\sigma_{x^{m}}$ and $I_{r}$ is the mean Cartan tensor of $F$.

Also, by the assumption, Randers metric $F$ is of isotropic $S$-curvature, that is, there is a scalar function $\lambda(x)$ on $M$ such that
\be
{\bf S}=(n+1)\lambda(x)F. \label{S2}
\ee
On the other hand,  by the assumption that Randers metric $\bar{F}:=e^{\sigma(x)}F$ is of weakly isotropic scalar curvature and by Theorem \ref{SCS}, $\bar{F}$ must be of isotropic $S$-curvature ,
\be
{\bar{\bf S}}=(n+1){\bar\lambda}(x){\bar F}, \label{S3}
\ee
where $\bar{\lambda}(x)$ is a scalar function on $M$. Plugging ($\ref{gij}$), ($\ref{Ii}$), ($\ref{S2}$) and ($\ref{S3}$) into ($\ref{S1}$), and then, similar to the proof of Theorem 1.2 in {\cite{CY}}, we can conclude that $\sigma(x)$ is a constant. \qed

\vskip 2mm

Lemma \ref{lmS} shows that,  if $F$ is a non-Riemannian Randers metric with isotropic $S$-curvature on an $n$-dimensional manifold $M (n\geq 3),$  then there is no non-constant scalar function $\sigma=\sigma(x)$ such that $\bar{F}:=e^{\sigma} F$ is of weakly isotropic scalar curvature.

Now, we are in the position to prove Theorem \ref{CWS}.

{\bf Proof of Theorem \ref{CWS}.} \ By the assumption, $F=\alpha+\beta$ is a conformally flat metric, that is, there exists  a scalar function $\kappa(x)$ on $M$, such that
\be
F=e^\kappa(x){\bar{F}},
\ee
where ${\bar{F}}$ is a Minkowski metric.

Obviously, $\bar{F}$ is of isotropic S-curvature, $\bar{\bf S}=0$. Further,  since $F=e^\kappa(x){\bar{F}}$ is of weakly isotropic scalar curvature,  by Lemma \ref{lmS}, $\kappa(x)$ must be a constant.  Hence, $F$ is a Minkowski metric. This completes the proof. \qed

\vskip 10mm

\noindent
Xinyue Cheng \\
School of Mathematical Sciences \\
Chongqing Normal University \\
Chongqing 401331,  P. R. China \\
chengxy@cqnu.edu.cn

\vskip 5mm

\noindent
Yannian Gong \\
School of Sciences \\
Chongqing University of Technology \\
Chongqing 400054,  P. R. China \\
gyn@2017.cqut.edu.cn

\end{document}